\documentclass[10pt]{article}
\usepackage{amsfonts,latexsym,amstext}
\usepackage{amsmath}
\usepackage{amssymb}
\usepackage[english]{babel}
\usepackage[latin1]{inputenc}

\usepackage[usenames]{color}
\definecolor{red}{rgb}{1.0,0.0,0.0}
\def\red#1{{\textcolor{red}{#1}}}
\definecolor{blu}{rgb}{0.0,0.0,1.0}

\definecolor{gre}{rgb}{0.03,0.50,0.03}

\newtheorem{theorem}{Theorem}[section]
\newtheorem{lemma}[theorem]{Lemma}
\newtheorem{proposition}[theorem]{Proposition}
\newtheorem{definition}[theorem]{Definition}

\newtheorem{hypothesis}[theorem]{Hypothesis}

\newtheorem{remark}[theorem]{Remark}
\newtheorem{corollary}[theorem]{Corollary}

\numberwithin{equation}{section}

\newcommand{\myref}[1]{(\ref {#1})}

\def\qed{{\hfill\hbox{\enspace${ \square}$}} \smallskip}
\def\sqr#1#2{{\vcenter{\vbox{\hrule height .#2pt \hbox{\vrule
 width .#2pt height#1pt \kern#1pt \vrule
width .#2pt} \hrule height .#2pt}}}}
\def\square{\mathchoice\sqr54\sqr54\sqr{4.1}3\sqr{3.5}3}

\def\eps{\varepsilon}

\def\ds{\begin{displaystyle}}
\def\eds{\end{displaystyle}}
\def\dis{\displaystyle }
\def\<{\left\langle }
\def\>{\right\rangle }

\def\dim{\noindent \hbox{{\bf Proof.} }}

\def\R{\mathbb R}
\def\N{\mathbb N}

\def\E{\mathbb E}
\def\P{\mathbb P}

\def\cald{{\cal D}}

\def\calf{{\cal F}}

\def\calh{{\cal H}}
\def\calk{{\cal K}}
\def\call{{\cal L}}
\def\caln{{\cal N}}
\def\calp{{\cal P}}

\def\calu{{\cal U}}

\def\call{{\cal L}}

\def\to{\rightarrow}

\begin{document}

\title{Stochastic Optimal Control with Delay in the Control II:\\
Verification Theorem and Optimal Feedbacks}
\date{}

\author{Fausto Gozzi
\\
Dipartimento di Economia e Finanza\\
Universit\`a LUISS - Guido Carli\\
Viale Romania 32,
00197 Roma,
Italy\\
e-mail: fgozzi@luiss.it\\
\\
Federica Masiero\\
Dipartimento di Matematica e Applicazioni\\ Universit\`a di Milano Bicocca\\
via Cozzi 55, 20125 Milano, Italy\\
e-mail: federica.masiero@unimib.it}

\maketitle
\begin{abstract}
We consider a stochastic optimal control problem governed by a stochastic differential equation with delay in the control.
Using a result
of existence and uniqueness of a sufficiently regular mild solution of the associated Hamilton-Jacobi-Bellman (HJB) equation, see
the companion paper \cite{FGFM},
we solve the control problem by proving a Verification Theorem and the existence of optimal feedback controls.
%
\end{abstract}

\bigskip \noindent


\bigskip \noindent

\textbf{Acknowledgements}:

Financial support from the grant
MIUR-PRIN 2010-11 ``Evolution differential problems: deterministic
and stochastic approaches and their interactions''
is gratefully acknowledged. The second author have been
supported by the Gruppo Nazionale per l'Analisi Matematica,
la Probabilit\`a e le loro Applicazioni (GNAMPA)
of the Istituto Nazionale di Alta Matematica (INdAM).

\newpage

\tableofcontents

\section{Introduction}
The aim of this paper is to solve, through the dynamic programming approach,
a class of stochastic optimal control problems with delay in the control.

Stochastic optimal control problems governed by delay equations with delay in the control are usually harder to study than the ones when the delay appears only in the state (see e.g. \cite[Chapter 4]{BDDM07} in the deterministic case and \cite{GM,GMSJOTA} in the stochastic case).
When one tries to apply the dynamic programming method the main difficulty is the fact that, even in the simplified setting introduced first by Vinter and Kwong \cite{VK} in the deterministic case (see e.g. \cite{GM} for the stochastic case), the associated HJB equation is an infinite dimensional second order semilinear PDE that
does not satisfy the so-called ``structure condition'', which substantially
means that the control can act on the system modifying its dynamics
at most along the same directions along which the noise acts.
The absence of such condition, together with the lack of smoothing properties which is a common feature of problems with delay,
prevents the use of the known techniques, based on BSDE's
(see e.g. \cite{FT2})
or on fixed point theorems in spaces of continuous functions
(see e.g. \cite{CDP1,CDP2,DP3,G1,G2})
or in Gauss-Sobolev spaces
(see e.g. \cite{ChowMenaldi,GGSPA}),
to prove the existence of regular solutions of this HJB equation.\footnote{The viscosity solution technique can still be used (see e.g. \cite{GMSJOTA})
but to prove existence (and possibly uniqueness) of solutions that are merely continuous.}

In the companion paper \cite{FGFM} we proved, using an ad hoc method based on the partial smoothing properties of the Ornstein-Uhlenbeck semigroup, that
the HJB equation associated to the class of problems under study here,
admits a unique {\em mild solution}, i.e. a solution which possesses enough regularity to give sense to the ``candidate optimal feedback map''
which depends on a suitable directional derivative (denoted by $\nabla^B$)
of the solution.

In this paper we start from such a result and exploit it to solve our class of problems. More precisely we prove here:
\begin{itemize}
  \item[(A)] an approximation result for the solutions of the HJB equation, i.e. that the mild solutions can be approximated by classical solutions, to which Ito's formula applies (Lemma \ref{lm:approximation});
  \item[(B)] a verification theorem for the control problem
  (Theorem \ref{teorema controllo});
  \item[(C)] the existence, under further assumptions on the Hamiltonian, of optimal controls in feedback form (Theorem \ref{teo su controllo feedback}).
\end{itemize}
These results allows to treat satisfactorily a specific class of state equations and data which arises naturally in many applied problems (see e.g. \cite{FabbriGozzi08,FedTacSicon,GM,GMSJOTA,Kol-Sha,PhamBruder}).

The three points outlined above are the ones followed, in a different context, in \cite{G1}
(see also Chapter 4 of \cite{FabbriGozziSwiech} and, in particular, Section 4.8). However the specific features of our problems prevents the use of the same techniques. In particular the approximation result, which is a key tool here,
must be completely reworked, as we explain in Remarks \ref{rm:defpistrictstrong}.
and \ref{rm:Yk}.

We finally note that, similarly to what often happens in the literature (see e.g.
\cite{FedTacSicon,GM,GMSJOTA}) here we treat the case of ``distributed delay''
which gives rise to a bounded control operator in the state equation.
The case of ``pointwise delay'', even if it seems treatable with our approach,
is left for future extensions of our research.

%
%


\subsection{Our results in a simple motivating case}
\label{SS:motiv}

To be more clear we now briefly describe our setting and our main result in a special case. Let $(\Omega, \calf,  \P)$
be a complete probability space
and consider the following linear controlled Stochastic
Differential Equation (SDE) in $\R$:
\begin{equation}
\left\{
\begin{array}
[c]{l}
dy(s)  =a_0 y(s) dt+b_0 u(s) ds +\dis\int_{-d}^0b_1(\xi)u(s+\xi)d\xi+\sigma dW_s
,\text{ \ \ \ }s\in[t,T] \\
y(t)  =y_0,\\
u(\xi)=u_0(\xi), \quad \xi \in [-d,0).
\end{array}
\right.  \label{eq-contr-ritINTRO}
\end{equation}
Here $W$ is a standard Brownian motion in $\R$, and $(\calf_t)_{t\geq 0}$ is the
augmented filtration generated by $W$.
We assume $a_0,b_0\in \R$, $\sigma >0$. The parameter $d>0$ represents the maximum delay
the control takes to affect the system while $b_1$ is the density function
taking account the aftereffect of the control on the system.
The case treated here is when $b_1\in L^2([-d,0],\R)$ (``distributed delay'')
while a more difficult case which we leave for further research is when $b_1$ is a measure, e.g. a Dirac delta in $-d$ (``pointwise delay'').

The initial data are the initial state $y_0$ and the past history $u_0$ of the control.
The control $u$ belongs to $L^2_{\calf}(\Omega\times [0,T], U)$,
the space of predictable square integrable processes with values in
$U\subseteq \R$, closed.

Such kind of equations is used e.g. to model the effect of advertising on
the sales of a product (see e.g. \cite{GM,GMSJOTA}), the effect of investments
on growth (see e.g. \cite{FabbriGozzi08} in a deterministic case),
or, in a more general setting, to model optimal portfolio problems with execution delay,
(see e.g. \cite{PhamBruder}) or to model the interaction of drugs with tumor cells
(see e.g. \cite{Kol-Sha} p.17 in the deterministic case).

In many applied cases (like the ones quoted above) the goal of the problem is to minimize the following objective functional
\begin{equation}\label{costoconcretoINTRO}
J(t,x_0,u_0;u(\cdot))=\E \int_t^T \left(\red{\bar\ell_0(s)}+\bar\ell_1(u(s))\right)\;ds +\E  \bar\phi(y(T)).
\end{equation}
where \red{$\bar\ell_0:[0,T]\rightarrow \R$} , $\bar\ell_1:U\rightarrow \R$ and
$\bar\phi:\R\rightarrow \R$ are continuous functions satisfying suitable assumptions that will be introduced in Section \ref{section-statement}.
It is important to note that here $\bar\ell_0$, $\bar\ell_1$ and $\bar\phi$ do not depend on the past of the state and/or control.
This is a very common feature of such applied problems.

A standard way\footnote{It must be noted that, under suitable restrictions on the data, one can treat optimal control
problems with delay in the control also by a direct approach without transforming them in infinite dimensional problems.
However in the stochastic case such direct approach seems limited to a very special class of cases
(see e.g. \cite{LarssenRisebro03}) which does not include our model and models commonly used in applications
like the ones just quoted.} to approach these delayed control problems, introduced in \cite{VK} for the deterministic case
(see \cite{GM} for the stochastic case) is to reformulate the above linear delay equation as a linear SDE in the Hilbert space
$\calh:=\R \times L^2([-d,0],\R)$, with state variable $Y=(Y_0,Y_1)$ as follows.
\begin{equation}
\left\{
\begin{array}
[c]{l}
dY(s)  =AY(s) ds+Bu(s) ds+GdW_s
,\text{ \ \ \ }s\in[t,T] \\
Y(t)=x=(x_0,x_1),
\end{array}
\right.   \label{eq-astrINTRO}
\end{equation}
where $A$ generates a $C_0$-semigroup (see next section for precise definitions)
while, at least formally,\footnote{Note that, when $b_1$ is a measure, the above operator
$B$ is not bounded and this makes the problem more difficult.}
\begin{equation}
 \label{BGINTRO}
B:\R\rightarrow \calh,\qquad Bu=(b_0 u, b_1(\cdot)u), \; u\in\R,
\qquad
G:\R\rightarrow \calh,\qquad Gx=(\sigma x, 0), \; x\in\R.
\end{equation}
Moreover $x_0=y_0$ while $x_1(\xi)=\int_{-d}^\xi b_1(\varsigma)u_0(\varsigma-\xi)d\varsigma$, ($\xi \in [-d,0]$) i.e.
the infinite dimensional datum $x_1$ depends on the past of the control.

The value function is defined as
$$
V(t,x):=\inf_{u(\cdot)\in L^2_{\calf}(\Omega\times [0,T], U)}
\E \left(\int_t^T \left[\red{\bar\ell_0(s)}+\bar\ell_1(u(s))\right] ds +\bar\phi(Y_0(T))
\right)
$$
The associated HJB equation (whose candidate solution is the value function) is
\begin{equation}\label{HJBINTRO}
\left\{\begin{array}{l}\dis
-\frac{\partial v(t,x)}{\partial t}=
\frac{1}{2}Tr \;GG^*\nabla^2v(t,x)
+ \< Ax,\nabla v(t,x)\>_\calh
+ \bar H_{min} (\nabla v(t,x)) +\red{\bar\ell_0(t)},\qquad t\in [0,T],\,
x\in \calh,\\
\\
\dis v(T,x)=\bar\phi(x_0),
\end{array}\right.
\end{equation}
where, defining the current value Hamiltonian $\bar H_{CV}$ as
$$
\bar H_{CV}(p;u):=\<p,Bu\>_\calh+\bar\ell_1(u)=\<B^*p,u\>_\R +\bar\ell_1(u)
$$
we have
\begin{equation}\label{HminINTRO}
\bar H_{min} (p):=
\inf_{u\in U} \bar H_{CV}(p\,;u).
\end{equation}

It is well known, see e.g. \cite{YongZhou99} Section 5.5.1, that, if the value function
$V$ is smooth enough, and if the current value Hamiltonian $\bar H_{CV}$ always admits at least a minimum point, a natural candidate
optimal feedback map is given by $(t,x) \mapsto u^*(t,x)$ where
$u^*(t,x)$ satisfies
\[
 \<\nabla V(t,x),Bu^*(t,x)\>_{\R}+\bar\ell_1(u^*(t,x))= \bar H_{min}(\nabla V(t,x)).
\]
i.e. where $u^*(t,x)$ is a minimum point of the function
$u\mapsto \bar H_{CV}(\nabla v(t,x);u)$, $\R\to \R$.
To take account of the presence of $B$ it is convenient to define,
for $z \in \R$,
\begin{equation}\label{HminINTRObis}
H_{min} (z):=
\inf_{u\in U}  \{\<z,u\>_\R+\bar\ell_1(u) \}=:\inf_{u\in U} H_{CV}(z\,;u) \}
\end{equation}
so that
\begin{equation}\label{HminINTROter}
\bar H_{min} (p)=
\inf_{u\in U} \bar H_{CV}(p\,;u).
=\inf_{u\in U} \{\<B^*p,u\>_\R+\bar\ell_1(u ) \}
=:H_{min}(B^*p).
\end{equation}
From now on we will use $H_{CV}$ and $H_{min}$ in place of $\bar H_{CV}$ and $\bar H_{min}$ writing $H_{min}(\nabla^B v(t,x))$ in place of $\bar H_{min}(\nabla v(t,x))$ in \myref{HJBINTRO}.

Since $u^*(t,x)$ is a minimum point of the function
$u\mapsto H_{CV}(B^*\nabla v(t,x);u)$, $\R\to \R$,
the minimal regularity required to give sense to such term is the existence
of $B^*\nabla v(t,x)$ which we will call $\nabla^B v(t,x)$ according to the definition and notation 
used e.g. in \cite{FTGgrad,Mas} and in the companion paper \cite{FGFM}.
In \cite{FGFM} a result of existence and uniqueness of mild solutions for
such equation has been proved. Mild solutions (defined through an integral form of
\myref{HJBINTRO} are continuous and such that
$\nabla^B v$ is well defined and continuous, hence the
natural candidate optimal feedback map $u^*(t,x)$ above is well defined.

Notice that if the controlled state equation satisfies the ``structure condition'', meaning that the control affects the system only through
the noise, then by the so called BSDEs approach, see e.g. \cite{FT2}, the fundamental relation and the consequent verification theorem,
can be proved also by applying the Girsanov Theorem. In our case the ``structure condition'' does not hold,
since it would mean that ${\rm Im} B\subseteq {\rm Im} G$. This is an intrinsic feature of control problems with delay in the control
since the fact that $\operatorname{Im} B $ is not contained in $\operatorname{Im} G$ is just due to the presence of the delay in the control.
If the delay in the control disappears, then the structure condition hold, even if delay in the state is present (see e.g. \cite{GM,GMSJOTA,FT2,MasBanach}).

Then solve the problem, as recalled in the beginning of this introduction,
one needs to accomplish the steps (A)-(B)-(C), that we briefly introduce:
\begin{itemize}
 \item [(A)] In Lemma \ref{lm:approximation} we show that if we suitably approximate the coefficients of the HJB equation (\ref{HJBINTRO}), we obtain
 a sequence of functions $(w_n)_n$ which are strict solutions of the approximating HJB equations, and which
 are once differentiable in time and twice differentiable with respect to $x$. Moreover the sequence $(w_n)_n$ converges to $v$, the mild solution
 of the HJB equation, in the sense of the $\calk$-convergence, see Definition \ref{k-conv}.
 \item [(B)] In Theorem \ref{teorema controllo} we apply the fundamental relation proved in Proposition \ref{prop rel fond}: in order to prove the fundamental relation
 it is crucial to apply the Ito formula, and this can be done thanks to the approximation performed in (A) and to a further approximation of the state.
  \item [(C)] In Theorem \ref{teo su controllo feedback}, under further regularity assumptions stated in Hypothesis
   \ref{ipotesicostoconcretobis}, we solve the closed loop equation and so we show the existence of optimal controls in feedback form and the fact that the value function coincides
   with the solution of the HJB equation, see Theorem \ref{teo:v=V}.
   \end{itemize}

\noindent

%

\bigskip

Finally notice that in the present paper we deal with a finite dimensional control delay equation
(\ref{eq-contr-ritINTRO}), that here in the introduction we have presented in dimension one for the sake of simplicity. The same arguments apply
if we consider the case of a controlled stochastic differential equation in an infinite dimensional Hilbert space $\calh_0$ with delay in the control
as follows.
\begin{equation}
\left\{
\begin{array}
[c]{l}%
dy(t)  =A_0 y(t) dt+B_0 u(t) dt +\int_{-d}^0B_1(\xi)u(t+\xi)d\xi+\sigma dW_t
,\text{ \ \ \ }t\in[0,T] \\
y(0)  =y_0,\\
u(\xi)=u_0(\xi), \quad \xi \in [-d,0).
\end{array}
\right.  \label{eq-contr-rit-infINTRO}
\end{equation}
Here $W$ is a cylindrical Wiener process in another Hilbert space $\Xi$, $A_0$ is the generator of a strongly continuous semigroup in
$\calh_0$, $\sigma \in \call (\Xi,\calh_0)$, and we have to assume some smoothing properties for the Ornstein Uhlenbeck transition
semigroup with drift term given by $A_0$ and diffusion equal to $\sigma$, see Remarks 3.2, 4.7 and 4.12 of \cite{FGFM} for more details.

\subsection{Plan of the paper}

In the first two Sections we present the control problem we treat and we recall some results on the state equation and on the HJB equation proved in \cite{FGFM},
in the last two Sections we solve the control problem. In details, the plan of the paper is the following:
\begin{itemize}
  \item in Section \ref{section-prel} we give some notations and we present the problem and the main assumptions;
  \item in Section \ref{section-smoothOU} we collect some results proved in \cite{FGFM} and fundamental in solving the control problem:
  the partial smoothing property for the Ornstein-Uhlenbeck transition semigroup related to the state equations, and results on
  the existence of a mild solution of the HJB equation;
\item in Section \ref{sec-Kstrong} we prove that this mild solution can be approximated by means of a sequence of strong solutions;
  \item Section \ref{sec-verifica} is devoted to solve the optimal control problem. In Subsection 5.1
  we prove a verification theorem, and finally in Subsection 5.2 we identify the value function of the control problem with the solution
  of the HJB equation and we characterize the optimal control by a feedback law.
\end{itemize}

\section{Some preliminary results on the control problem}\label{section-prel}

In this section we collect and synthesize, for the reader's convenience, the basic material on our control problem which has been already exposed in the companion paper \cite{FGFM}.

\subsection{Notation and $C$-derivatives}\label{subsection-notation}

Let $H$, $K$ be real and separable Hilbert spaces, with norms given by $\left|  x\right|_{H}$
and by $\left|  x\right|_{K}$, respectively, or by $\left|  x\right|  $, if no
confusion is possible, and by scalar product $\left\langle \cdot,\cdot\right\rangle _H$,
$\left\langle \cdot,\cdot\right\rangle _K$, respectively, or simply
by $\left\langle \cdot,\cdot\right\rangle $. The space $\call(H,K)$ denotes the
space of bounded linear operators from $H$ to $K$ endowed with the usual
operator norm.


In the following, by $(\Omega, \calf, \P)$ we denote a complete probability
space, and by $L^2_\calp(\Omega\times[0,T],H)$
the Hilbert space of all predictable processes
$(Z_t)_{t\in[0,T]}$ with values in $H$, normed by
$\Vert Z\Vert^2 _{L^2_\calp(\Omega\times[0,T],H)}=\E\int_0^T\vert Z_t\vert^2\,dt$.

Next we introduce some spaces of functions. We let $H$ and $Z$ be Hilbert spaces.
By $B_b(H,Z)$ (respectively $C_b(H,Z)$, $UC_b(H,Z)$) we denote the space of all functions
$f:H\rightarrow Z$ which are Borel measurable and bounded (respectively continuous
and bounded, uniformly continuous and bounded).

Given an interval $I\subseteq \R$ we denote by
$C(I\times H,Z)$ (respectively $C_b(I\times H,Z)$)
the space of all functions $f:I \times H\rightarrow Z$
which are continuous (respectively continuous and bounded).
$C^{0,1}(I\times H,Z)$ is the space of functions
$ f\in C(I\times H)$ such that for all $t\in I$
$f(t,\cdot)$ is Fr\'echet differentiable.
By $UC_{b}^{1,2}(I\times H,Z)$
we denote the linear space of the mappings $f:I\times H \to Z$
which are uniformly continuous and bounded
together with their first time derivative $f_t$ and its first and second space
derivatives $\nabla f,\nabla^2f$.
If $Z=\R$ we omit it in all the above spaces.

We will need spaces of functions which are differentiable in suitable directions, we refer to the definition of $C$-directional derivatives
given in \cite{Mas}, Section 2, in \cite{FTGgrad} and in \cite{FGFM}, Section 2. More precisely, we refer to \cite{FGFM}, Definition 2.1,
for the definiton of of
$C$-directional derivative, and of $C$-G\^ateaux or $C$-Fr\'echet differentiable function; we refer to \cite{FGFM}, Definition 2.2, for the definitions of the
spaces of functions $C^{1,C}_{b}(H,Z)$,
$C^{0,1,C}_b(I\times H,Z)$,  $C^{0,2,C}_b(I\times H,Z)$.
\newline We also recall that, according to \cite{FGFM}, Definition 2.2,
for any $\alpha\in(0,1)$ and $T>0$ (this time $I$ is equal to $[0,T]$) we denote
by $C^{0,1}_{\alpha}([0,T]\times H)$ the space of functions
$ f\in C_b([0,T]\times H))\cap C^{0,1}((0,T]\times H,Z)$ such that
the map $(t,x)\mapsto t^{\alpha} \nabla f(t,x)$
belongs to $C_b((0,T]\times H,H)$ and by
 We also denote
by $C^{0,1,C}_{\alpha}([0,T]\times H)$ the space of functions
$ f\in C_b([0,T]\times H)\cap C^{0,1,C}((0,T]\times H)$ such that
the map $(t,x)\mapsto t^{\alpha} \nabla^C f(t,x)$
belongs to $C_b((0,T]\times H,K)$.
\newline Again by \cite{FGFM}, Definition 2.2,  for any $\alpha\in(0,1)$ and $T>0$)
we denote  by $C^{0,2}_{\alpha}([0,T]\times H)$
the space of functions
$ f\in C_b([0,T]\times H)\cap C^{0,2}((0,T]\times H)$ such that 
the map $[0,T]\times H \to H$,
$(t,x)\mapsto \nabla f(t,x)$ is bounded and continuous
and the map $(0,T]\times H \to \call( H,H)$,
$(t,x)\mapsto t^{\alpha} \nabla^2 f(t,x)$ is bounded and continuous,
and we denote by $C^{0,2,C}_{\alpha}([0,T]\times H)$
the space of functions
$ f\in C_b([0,T]\times H)\cap C^{0,2,C}((0,T]\times H)$ such that 
{the map $[0,T]\times H \to H$,
$(t,x)\mapsto \nabla f(t,x)$ is bounded and continuous
and the map $(0,T]\times H \to \call( K,H)$,
$(t,x)\mapsto t^{\alpha} \nabla^C\nabla f(t,x)$ is bounded and continuous.
\newline The spaces $C^{0,1}_{\alpha}([0,T]\times H),\,C^{0,1,C}_{\alpha}([0,T]\times H),\,C^{0,2}_{\alpha}([0,T]\times H)$ and
$C^{0,2}_{\alpha}([0,T]\times H)$ are Banach spaces if endowed with suitble norms, see \cite{FGFM}, Definition 2.2.

\subsection{Setting of the problem and main assumptions}\label{section-statement}

In this paper we are concerned with the solution of an optimal control problem related to
an $n$-dimensional controlled equation with delay in the control, that we are going to introduce.\\
In a complete probability space $(\Omega, \calf,  \P)$
we consider
\begin{equation}
\left\{
\begin{array}
[c]{l}%
dy(t)  =a_0 y(t) dt+b_0 u(t) dt +\dis\int_{-d}^0b_1(\xi)u(t+\xi)d\xi+\sigma dW_t
,\text{ \ \ \ }t\in[0,T] \\
y(0)  =y_0,\\
u(\xi)=u_0(\xi), \quad \xi \in [-d,0),
\end{array}
\right.  \label{eq-contr-rit}
\end{equation}
where we assume the following.

\begin{hypothesis}\label{ipotesibasic}
\begin{itemize}
\item[]
  \item[(i)] $W$ is a standard Brownian motion in $\R^k$, and $(\calf_t)_{t\geq 0}$ is the
augmented filtration generated by $W$;
  \item[(ii)] $a_0\in \call(\R^n;\R^n)$, $\sigma$ is in $\call(\R^k;\R^n)$;
  \item[(iii)] the control strategy $u$ belongs to $\calu$ where
$$\calu:=\left\lbrace z\in L^2_{\calp}(\Omega\times [0,T], \R^m):
u(t)\in U \;a.s.\right\rbrace $$
where $U$ is a closed subset of $\R^n$;
  \item[(iv)] $d>0$ (the maximum delay the control takes to affect the system);
  \item[(v)] $b_0 \in \call(\R^m;\R^n)$;
  \item[(vi)]
$b_1\in L^2([-d,0],\call(\R^m;\R^n)).$
\end{itemize}
\end{hypothesis}

\medskip

\subsection{Infinite dimensional reformulation}

Following the approach of \cite{VK}, applied in \cite{GM} to the stochastic case, we reformulate equation (\ref{eq-contr-rit}) as an abstract stochastic
differential equation in the Hilbert space $\calh=\R^n\times L^2([-d,0],\R^n)$.
To this end we introduce the operator $A : \cald(A) \subset \calh
\rightarrow \calh$ as follows: for $(y_0,y_1)\in \calh$
\begin{equation}\label{A}
A(y_0 ,y_1 )=( a_0 y_0 +y_1(0), -y_1'), \qquad
\cald(A)=\left\lbrace(y_0,y_1)\in \calh:y_1\in W^{1,2}([-d,0],\R^n), y_1(-d)=0 \right\rbrace.
\end{equation}
We denote by $e^{tA}$ the $C_0$-semigroup generated by $A$: for
$y=(y_0,y_1)\in \calh$,
\begin{equation}
e^{tA} \left(\begin{array}{l}y_0 \\y_1\end{array}\right)=
\left(
\begin{array}
[c]{ll}%
e^{ta_0 }y_0+\int_{-d}^{0}1_{[-t,0]} e^{(t+s)a_0 } y_1(s)ds \\[3mm]
y_1(\cdot-t)1_{[-d+t,0]}(\cdot).
\end{array}
\right)  \label{semigroup}
 \end{equation}
  We will use, for $N\in \N$ big enough, the resolvent operator
 $(N-A)^{-1}$ which can be computed explicitly
 \begin{equation}
 (N-A)^{-1} \left(\begin{array}{l}y_0 \\y_1\end{array}\right)=
 \left(
 \begin{array}
 [c]{ll}%
 (N-a_0)^{-1}\left[y_0+\int_{-d}^{0} e^{N s } y_1(s)ds\right]
 \\[3mm]
 \int_{-d}^{\cdot} e^{N(s-\cdot)} y_1(s)ds.
 \end{array}
 \right)  \label{resolvent}
 \end{equation}
Similarly, denoting by $e^{tA^*}=(e^{tA})^*$ the $C_0$-semigroup generated by $A^*$,
we have for
$z=\left(z_0,z_1\right)\in \calh $
\begin{equation}
e^{tA^*} \left(\begin{array}{l}z_0 \\z_1\end{array}\right)=
\left(
\begin{array}[c]{ll}
e^{t a_0^* }z_0 \\[3mm]
e^{(\cdot+t) a_0^* }z_0 1_{[-t,0]}(\cdot) +z_1(\cdot+t)1_{[-d,-t)}(\cdot).
\end{array}
\right)  \label{semigroupadjoint}
\end{equation}
The infinite dimensional noise operator is defined as
\begin{equation}
 \label{G}
G:\R^{k}\rightarrow \calh,\qquad Gy=(\sigma y, 0), \; y\in\R^k.
\end{equation}
The control operator $B$
is bounded and defined as
\begin{equation}
 \label{B}
B:\R^{m}\rightarrow \calh,\qquad Bu=(b_0 u, b_1(\cdot)u), \; u\in\R^m
\end{equation}
and its adjoint is
\begin{equation}
 \label{B*}
B^*:\calh^* \rightarrow \R^{m},\qquad B^*(x_0,x_1)=
b^*_0 x_0+\int_{-d}^0 b_1(\xi)^*x_1(\xi)d\xi, \; (x_0,x_1)\in\calh.
\end{equation}
Note that, in the case of pointwise delay the last term of the drift in
the state equation \myref{eq-contr-rit} is $u(t-d)$, hence $b_1(\cdot)$ is
a measure: the Dirac delta $\delta_{-d}$.
Hence in this case $B$ is unbounded as it takes values in
$\R^n \times C^*([-d,0],\R^n)$ (here we denote by $C^*([-d,0],\R^n)$ the dual space of $C([-d,0],\R^n)$).
%
%
%

Given any initial datum $(y_0,u_0)\in \calh$ and any admissible control $u\in \calu$ we call $y(t;y_0,u_0,u)$ (or simply $y(t)$ when clear from the context) the unique solution (which comes from standard results on SDE's, see e.g.  \cite{IkedaWatanabe} Chapter 4, Sections 2 and 3)
of (\ref{eq-contr-rit}).

Let us now define the process
$Y=(Y_0,Y_1)\in L^2_\calp(\Omega \times [0,T],\calh)$ as
$$
Y_0(t)=y(t), \qquad Y_1(t)(\xi)=\int_{-d}^\xi u(\zeta+t-\xi)b_1(\zeta)d\zeta,
$$
where $y$ is the solution of equation (\ref{eq-contr-rit}) and $u$ is the control process in (\ref{eq-contr-rit}).
By Proposition 2 of \cite{GM},
the process $Y$
is the unique solution of the abstract evolution equation
in $\calh$
\begin{equation}
\left\{
\begin{array}
[c]{l}
dY(t)  =AY(t) dt+Bu(t) dt+GdW_t
,\text{ \ \ \ }t\in[ 0,T] \\
Y(0)  =y=(y_0,y_1),
\end{array}
\right.   \label{eq-astr}%
\end{equation}
where $y_0=x_0$ and $y_1(\xi)=\dis\int_{-d}^\xi u_0(\zeta-\xi)b_1(\zeta)d\zeta$.
Note that we have $y_1\in L^2([-d,0];\R^n)$
Taking the integral (or mild) form of (\ref{eq-astr}) we have
\begin{equation}
Y(t)  =e^{tA}y+\int_0^te^{(t-s)A}B u(s) ds +\int_0^te^{(t-s)A}GdW_s
,\text{ \ \ \ }t\in[ 0,T]. \\
  \label{eq-astr-mild}%
\end{equation}

\subsection{Optimal Control problem}


The objective is to minimize, over all controls in $\calu$,
the following finite horizon cost:
 \begin{equation}\label{costoconcreto}
J(t,x;u)=\E \int_t^T \left(\bar\ell_0(s,y(s))+\bar\ell_1(u(s))\right)\;ds +\E  \bar\phi(x(T)).
\end{equation}
where {$\bar\ell_0:[0,T]\rightarrow \R$ }and
$\bar\phi:\R^n\rightarrow \R$ are continuous and bounded while
$\bar\ell_1:U\rightarrow\R$ is measurable and bounded from below.
Referring to the abstract formulation (\ref{eq-astr}) the cost in (\ref{costoconcreto}) can be rewritten also as
\begin{equation}\label{costoconcreto1}
J(t,x;u)=\E \left(\int_t^T \left[\ell_0(s,Y(s))+\ell_1(u(s))\right]\,ds + \phi(Y(T))\right),
\end{equation}
where
$\ell_0:[0,T]\times \calh\rightarrow \R$, $\ell_1:U\to \R$ are defined by setting
\begin{equation}\label{l_0}
\red{\ell_0(t):=\bar\ell_0(t)} ,
\end{equation}
\begin{equation}\label{l_1}
\ell_1:=\bar\ell_1
\end{equation}
(here we cut the bar only to keep the notation homogeneous)
while $\phi :\calh\rightarrow \R$ is defined as
\begin{equation}\label{fi0}
\phi(x):=\bar\phi(x_0) \quad
\forall x=(x_0,x_1)\in \calh.
\end{equation}

Clearly, under the assumption above,
$\ell_0$ and $\phi$
are continuous and bounded while $\ell_1$ is measurable and bounded from below.
The value function of the problem is
\begin{equation}\label{valuefunction}
 V(t,x):= \inf_{u \in \calu}J(t,x;u).
\end{equation}
As done in Subsection \ref{SS:motiv}, we define the Hamiltonian in a modified way
(see \myref{HminINTRObis}); indeed, for $p\in \calh$, $u \in U$,
we define the current value Hamiltonian $H_{CV}$ as
$$
H_{CV}(p\,;u):=\<p,u\>_{\R^m}+\ell_1(u)
$$
and the (minimum value) Hamiltonian by
\begin{equation}\label{psi1}
H_{min}(p)=\inf_{u\in U}H_{CV}(p\,;u),
 \end{equation}
The associated HJB equation with unknown $v$ is then formally written as
\begin{equation}\label{HJBformale1}
  \left\{\begin{array}{l}\dis
-\frac{\partial v(t,x)}{\partial t}=\frac{1}{2}Tr \;GG^*\nabla^2v(t,x)
+ \< Ax,\nabla v(t,x)\>_\calh +\ell_0(t,x)+ {H}_{min} (\nabla^Bv(t,x)),\qquad t\in [0,T],\,
x\in D(A),\\
\\
\dis v(T,x)=\phi(x).
\end{array}\right.
\end{equation}

Existence of mild solutions of \myref{HJBformale1} is proved in \cite{FGFM}, and the following assumptions are needed.
\begin{hypothesis}\label{ipotesicostoconcreto}
\begin{itemize}
\item[]
  \item[(i)] $\phi\in C_b(\calh)$ and it is given by \myref{fi0} for
  a suitable $\bar\phi \in  C_b(\R^n)$;
  \item[(ii)] \red{$\ell_0\in C_b([0,T] )$} and it is given by \myref{l_0} ;
  \item[(iii)] $\bar\ell_1:U\rightarrow\R$ is measurable and bounded from below;
  \item[(iv)] the Hamiltonian $H_{min}:\R^m \to \R$ is Lipschitz continuous so
  there exists $L>0$ such that
  \begin{equation}\label{eq:Hlip}
    \begin{array}{c}
\vert H_{min }(p_1)-H_{min }(p_2)\vert\leq L \vert p_1-p_2\vert
\quad \forall\,p_1,\,p_2\in\R^m;
       \\[1.5mm]
\vert H_{min }(p)\vert\leq L(1 + \vert p\vert )
\quad \forall\,p\in\R^m.
    \end{array}
\end{equation}
\end{itemize}
\end{hypothesis}
To get more regular solutions (well defined second derivative $\nabla^B\nabla$,
which will be used to prove existence of optimal feedback controls) we will need
the following further assumption.

\begin{hypothesis}\label{ipotesicostoconcretobis}
  The Hamiltonian $H_{min}:\R^m \to \R$ is continuously differentiable
  and, for a given $L>0$, we have, beyond \myref{eq:Hlip},
  \begin{equation}\label{eq:Hlipder}
    \begin{array}{c}
\vert \nabla H_{min }(p_1)-\nabla H_{min }(p_2)\vert\leq L \vert p_1-p_2\vert
\quad \forall\,p_1,\,p_2\in\R^m;
    \end{array}
\end{equation}
\end{hypothesis}

\section{The Ornstein-Uhlenbeck semigroup and the HJB equation}
\label{section-smoothOU}
In the setting of Section \ref{section-statement}
we assume that Hypothesis \ref{ipotesibasic} holds true.
We take $\calh=\R^n \times L^2(-d,0;\R^n)$, $\Xi=\R^k$,
$(\Omega, \calf,  \P)$ a complete probability space, $W$ a standard
Wiener process in $\Xi$, $A$ and $G$ as in (\ref{A}) and (\ref{G}).
Then, for $x\in \calh$, we take the Ornstein-Uhlenbeck process $X^x(\cdot)$
given by

\begin{equation}
\left\{
\begin{array}
[c]{l}%
dX(t)  =AX(t) dt+GdW_t
,\text{ \ \ \ }t\ge 0\\
X(0)  =x,
\end{array}
\right.\label{ornstein-gen}
\end{equation}
In mild form, the Ornstein-Uhlenbeck process $X^x$ is given by
\begin{equation}
X^x(t)  =e^{tA}x +\int_0^te^{(t-s)A}GdW_s
,\text{ \ \ \ }t\ge 0. \\
  \label{ornstein-mild-gen}
\end{equation}
$X$ is a Gaussian process, namely for every $t>0$, the law of
$X(t)$ is $\caln (e^{tA}x,Q_t)$, the Gaussian measure with mean $e^{tA}x$ and
covariance operator $Q_t$,
where
\[
 Q_t=\int_0^t e^{sA}GG^*e^{sA^*}ds.
\]
The associated Ornstein-Uhlenbeck transition semigroup $R_t$ is defined by setting, for all $f\in B_b(\calh)$,
\begin{equation}
 \label{ornstein-sem-gen}
R_t[f](x)=\E f(X^x(t))
=\int_K f(z+e^{tA}x)\caln(0,Q_t)(dz).
\end{equation}
Given any $\bar\phi\in B_b(\R^n)$, we define, as in (\ref{fi0})
a function $\phi \in B_b(\calh)$, by setting
\begin{equation}\label{fi}
\phi(x)=\bar\phi(x_0) \quad
\forall x=\left(
x_0 ,x_1  \right)\in \calh.
\end{equation}
For such functions, the Ornstein-Uhlenbeck semigroup $R_t$ is written as
\begin{equation}
 \label{ornstein-sem-spec}
R_t[\phi](x)=\E \phi (X^x(t))=\E \bar\phi ((X^x(t))_0)=\int_\calh\bar\phi((z+e^{tA}x)_0)\caln(0,Q_t)(dz).
\end{equation}
For the Ornstein Uhlenbeck transition semigroup we have the following regularizing property, see also \cite{FGFM}, Propositions 4.9 and 4.11
for more details.

\begin{proposition}\label{lemmaderhpdeb}
Assume that Hypothesis \ref{ipotesibasic} holds. Assume moreover that, either
\begin{equation}\label{eq:hpdebreg}
\operatorname{Im}(e^{ta_0}b_0)\subseteq
\operatorname{Im}\sigma, \; \forall t > 0;
\qquad
\operatorname{Im}b_1(s)\in\operatorname{Im}\sigma,
\quad  a.e.\, \forall s\in[-d,0]
\end{equation}
or
\begin{equation}\label{eq:hpdebregbis}
\operatorname{Im}\left(e^{ta_0}b_0 +\int_{-d}^0 1_{[-t,0]}e^{(t+r)a_0}b_1(dr)
\right)
\subseteq\operatorname{Im}\sigma,
\quad \forall t>0.
\end{equation}
Then, for any bounded measurable $\phi$ as in (\ref{fi}),
$R_{t}\left[\phi\right]$ is $B$-Fr\'echet differentiable for every
$t>0$,
and,
 for every $k\in \R^m$,
for all $T>0$ there exists $C_T$ such that
\begin{equation}
\vert\<  \nabla^{B}(R_{t}\left[\phi\right])  (x),k\>_{\R^m}\vert
\leq C_T t^{-{1/2}} \Vert \bar\phi \Vert_\infty \;\vert k\vert_{\R^m}.
\label{eq:stimaderBnew}
\end{equation}
\end{proposition}

\subsection{Regular mild solutions of the HJB equation}
\label{sec-HJB}
We now recall results proved in \cite{FGFM} about existence of mild solutions of the HJB equation (\ref{HJBINTRO}).
We also state results about the existence of second order derivatives that will be needed in Section
\ref{sec:contr-feedback} to solve our control problem.

 First of all we introduce some suitable spaces of differentiable functions in $C_b([0,T]\times\calh)$. We fix $\alpha \in (0,1)$, in the following
 we will consider these spaces with $\alpha=1/2$.
\begin{definition}\label{df:Sigma}
Let $T>0$, $\alpha \in (0,1)$.
A function $g\in C_b([0,T]\times \calh)$ belongs to $\Sigma^1_{T,\alpha}$ if there exists a function$f\in C_{b}([0,T]\times \R^n)$
such that
$$g(t,x)=f\left(t,(e^{tA}x)_0\right),
\qquad \forall (t,x) \in [0,T]\times \calh\red{,}
$$
if, for any $t\in(0,T]$ the function $g(t,\cdot)$ is
$B$-Fr\'echet differentiable and if it exists
a function $\bar f\in C_b((0,T]\times \R^n;\R^m)$
such that
$$
t^\alpha \nabla^B g(t,x)=\bar f\left(t,(e^{tA}x)_0\right),
\qquad \forall (t,x) \in [0,T]\times \calh.
$$
\end{definition}
If, in the above definition, we take $f\in C_{\alpha}^{0,1}([0,T]\times \R^n)$,
then for any $t\in(0,T]$ the function $g(t,\cdot)$ is
both Fr\'echet differentiable and $B$-Fr\'echet differentiable.
Moreover, for $(t,x)\in [0,T]\times \calh$, $h \in \calh$, $k\in \R^m$,
$$
\<\nabla g(t,x),h\>_{\calh}=\<\nabla f\left(t,(e^{tA}x)_0\right),(e^{tA}h)_0\>_{\R^n},
\quad and \quad
\<\nabla^B g(t,x),k\>_{\R^m}=\<\nabla f\left(t,(e^{tA}x)_0\right),(e^{tA}Bk)_0\>_{\R^n}.
$$
This in particular implies that, for all $k\in \R^m$
\begin{equation}
\label{eq:nablaperSigma}
\<\nabla^B g(t,x),k\>_{\R^m}=\<\nabla g(t,x),Bk\>_\calh
\end{equation}
which also means $B^*\nabla g = \nabla^B g$.
Moreover, in such a case, the function $\bar f\in C_b((0,T]\times \R^n;\R^m)$
is given by
$$
\<\bar f(t,y),k\>_{\R^m}=
t^\alpha\<\nabla f\left(t,y\right),(e^{tA}Bk)_0\>_{\R^n},
\qquad (t,y)\in (0,T]\times \R^n, \quad k \in \R^m.
$$
What is written above still works if $B$ is unbounded:
it suffices to require $(e^{tA}B)_0$
bounded and continuous.
The set $\Sigma^1_{T,\alpha}$ is a closed subspace of $C^{0,1,B}_{\alpha}([0,T]\times \calh)$.

Next, in analogy to what we have done defining $\Sigma^1_{T,\alpha}$,
we introduce a subspace $\Sigma^2_{T,\alpha}$ of functions
$g \in C^{0,2,B}_{\alpha}([0,T]\times \calh)$
that depends in a special way on the variable $x\in\calh$.
\begin{definition}\label{df:Sigma2}
Let $T>0$, $\alpha \in (0,1)$.
A function $g\in C_b([0,T]\times \calh)$
belongs to $\Sigma^2_{T,\alpha}$ if:
\begin{itemize}
\item there exists a function
$f\in C_{b}([0,T]\times \R^n)$ such that
$$g(t,x)=f\left(t,(e^{tA}x)_0\right),
\qquad \forall (t,x) \in [0,T]\times \calh;
$$
\item for any $t\in [0,T]$, the function $g(t,\cdot)$ is
Fr\'echet differentiable and it exists
a function $\bar f\in C_b([0,T]\times \R^n;\calh)$ such that
$$
\nabla g(t,x)=\bar f\left(t,(e^{tA}x)_0\right),
\qquad \forall (t,x) \in [0,T]\times \calh;
$$
\item for any $t\in(0,T]$ the function $\nabla g(t,\cdot)$ is
$B$-Fr\'echet differentiable and it exists
a function $\bar{\bar f}\in C_b((0,T]\times \R^n;\call(\R^m,\calh))$
such that
$$
t^\alpha \nabla^B\nabla g(t,x)=\bar{\bar f}\left(t,(e^{tA}x)_0\right),
\qquad \forall (t,x) \in (0,T]\times \calh.
$$
\end{itemize}
\end{definition}
If, in the above definition, we take $f \in C^{0,2}_\alpha([0,T]\times \R^n)$ then, for any $t\in [0,T]$,
$$
\<\nabla g(t,x),h\>_\calh=
\<\nabla f\left(t,(e^{tA}x)_0\right),(e^{tA}h)_0\>_{\R^n},
\quad \hbox{for $(t,x)\in [0,T]\times \calh$, $h \in \calh$.}
$$
Moreover also $\nabla g(t,\cdot)$ is 
both Fr\'echet differentiable and $B$-Fr\'echet differentiable and
$$
\<\nabla^B\left(\<\nabla g(t,x),h\>\right),k\>_{\R^m}=
\<\nabla^2 f\left(t,(e^{tA}x)_0\right)(e^{tA}h)_0,(e^{tA}Bk)_0 \>_{\R^n},
\quad \hbox{for $(t,x)\in [0,T]\times \calh$, $h \in \calh$, $k \in \R^m$.}
$$
We also notice that, since the function $f$ is twice continuously
Fr\'echet differentiable the second order derivatives
$\nabla^B\nabla g$ and $\nabla\nabla^B g$ both exist and coincide:
$$
 \<\nabla^B\<\nabla g(t,x),h\>_{\calh},k\>_{\R^m}
 =\<\nabla\<\nabla^B g(t,x),k\>_{\R^m},h\>_{\calh}.
$$
Moreover, in such a case, the function 
$\bar f\in C_b([0,T]\times \R^n;\R^m)$
is given by
$$
\<\bar f(t,y),h\>_{\calh}=
\<\nabla f\left(t,y\right),(e^{tA}Bh)_0\>_{\R^n},
\qquad (t,y)\in [0,T]\times \R^n, \quad h \in\calh.
$$
Similarly, the function 
$\bar{\bar f}\in C_b\left((0,T]\times \R^n;\call(\R^m,\calh)\right)$
is given by
$$
\<\bar{\bar f}(t,y)k,h\>_{\calh}=
t^\alpha\<\nabla^2 f\left(t,y\right)(e^{tA}h)_0,(e^{tA}Bk)_0 \>_{\R^n}
\qquad (t,y)\in [0,T]\times \R^n, \quad h \in \calh, \; k\in R.
$$

When (\ref{eq:hpdebreg}) or
(\ref{eq:hpdebregbis}) hold we know, by Proposition \ref{lemmaderhpdeb},
that the function $g(t,x)=R_t[\phi](x)$
for $\phi$ given by (\ref{fi}) with $\bar\phi$ bounded and continuous,
belongs to $\Sigma^1_{T,1/2}$.
Now we introduce mild solutions to the HJB equation (\ref{HJBformale1}). By applying formally the variation of
constants formula to (\ref{HJBformale1}) we get the integral equation satisfied by the mild solution, that we rewrite here for the reader
convenience:
\begin{equation}
v(t,x) =R_{T-t}[\phi](x)+\int_t^T \left[R_{s-t}\left[
H_{min}(\nabla^B v(s,\cdot))\right](x)+\red{\ell_0(s)}\right]\; ds,\qquad t\in [0,T],\
x\in H.\label{solmildHJB}
\end{equation}
We use this formula to give the notion of mild solution for the HJB equation (\ref{HJBformale1}).

\begin{definition}\label{defsolmildHJB}
We say that a
function $v:[0,T]\times \calh\rightarrow\mathbb{R}$ is a mild
solution of the HJB equation (\ref{HJBformale1}) if the following
are satisfied:
\begin{enumerate}
\item $v(T-\cdot, \cdot)\in C^{0,1,B}_{{1/2}}\left([0,T]\times \calh\right)$;

\item  equality (\ref{solmildHJB}) holds on $[0,T]\times \calh$.
\end{enumerate}
\end{definition}

\begin{theorem}\label{esistenzaHJB}
Let Hypotheses \ref{ipotesibasic} and \ref{ipotesicostoconcreto} hold and
let (\ref{eq:hpdebreg}) or (\ref{eq:hpdebregbis}) hold.
Then the HJB equation (\ref{HJBformale1})
admits a mild solution $v$ according to Definition \ref{defsolmildHJB}.
Moreover $v$
is unique among the functions $w$ such that $w(T-\cdot,\cdot)\in\Sigma_{T,1/2}$ and it satisfies, for suitable $C_T>0$, the estimate
\begin{equation}\label{eq:stimavmainteo}
\Vert v(T-\cdot,\cdot)\Vert_{C^{0,1,B}_{{1/2}}}\le C_T\left(\Vert\bar\phi \Vert_\infty
+\Vert\bar\ell_0 \Vert_\infty \right).
\end{equation}
If the initial datum $\phi$ is also continuously $B$-Fr\'echet
(or Fr\'echet) differentiable,
then $v \in C^{0,1,B}_{b}([0,T]\times \calh)$ and, for suitable $C_T>0$,
\begin{equation}\label{eq:stimavmainteobis}
\Vert v\Vert_{C^{0,1,B}_{b}}\le C_T\left(\Vert\phi \Vert_\infty
+\Vert\nabla^B\phi \Vert_\infty+\Vert\ell_0 \Vert_\infty \right)
\end{equation}
(substituting $\nabla^B\phi$ with $\nabla\phi$ if $\phi$ is Fr\'echet differentiable).
If $\sigma$ is onto, then $v(T-\cdot, \cdot)\in
 C_{1/2}^{0,1}([0,T]\times\calh)$
\begin{equation}\label{stimadiffle-solHJB}
 \left\Vert v(T-\cdot,\cdot)\right\Vert _{C^{0,1}_{{1/2}}
  }\leq C_T\left(\Vert\bar\phi \Vert_\infty
+\Vert\bar\ell_0 \Vert_\infty \right)
\end{equation}
If also Hypothesis \ref{ipotesicostoconcretobis} holds then we have the following
\begin{itemize}
  \item[(i)] If $\bar\phi\in C^1_b(\R^n)$,
then we have $v \in \Sigma^2_{T,{1/2}}$,
hence the second order derivatives $\nabla^B\nabla v$ and $\nabla\nabla^B v$
exist and are equal.
Moreover there exists a constant $C>0$ such that\red{
\begin{equation}\label{stimanablav}
\vert \nabla v(t,x)\vert \leq C\Vert \nabla\bar\phi\Vert_\infty
\end{equation}
\begin{equation}\label{stimanablav^2}
\vert \nabla^B\nabla v(t,x)\vert= \vert \nabla\nabla^B v(t,x)\vert
\leq C (T-t)^{-{1/2}}\Vert \nabla\bar\phi\Vert_\infty
\end{equation}}
Finally, if $\sigma$ is onto, then also $\nabla^2 v$ exists and is continuous
and, for suitable $C>0$,\red{
\begin{equation}\label{stimanablav^2bis}
\vert \nabla^2 v(t,x)\vert
\leq C (T-t)^{-{1/2}}
\Vert \nabla\bar\phi\Vert_\infty
\end{equation}
}
  \item[(ii)]If only $\bar\phi \in C_b(\R^n)$
  then the function $(t,x)\mapsto (T-t)^{1/2}\nabla^B v(t,x)$ belongs to $\Sigma^1_{T,{1/2}}$
I
Moreover there exists a constant $C>0$ such that \red{
\begin{equation}\label{stimanablavreg}
\vert \nabla^B v(t,x)\vert \leq C(T-t)^{-1/2} \Vert\bar\phi\Vert_\infty
\end{equation}
\begin{equation}\label{stimanablav^2reg}
\vert \nabla^B\nabla^B v(t,x)\vert
\leq C (T-t)^{-1}\Vert \bar\phi\Vert_\infty
\end{equation}}
Finally, if $\sigma$ is onto, then also $\nabla^2 v$ exists and is continuous
in $[0,T)\times \calh$
and, for suitable $C>0$,\red{
\begin{equation}\label{stimanablav^2bisreg}
\vert \nabla^2 v(t,x)\vert
\leq C (T-t)^{-1}
\Vert \nabla\bar\phi\Vert_\infty
\end{equation}}
\end{itemize}
\end{theorem}

\section{$\calk$-convergence and approximations of solutions of the HJB equation}
\label{sec-Kstrong}

We first introduce the notion of $\calk$-convergence, following \cite{Ce95} and \cite{G1}
(see also Section B.5.1 of \cite{FabbriGozziSwiech}).
\begin{definition}\label{k-conv}
A sequence $(f_n)_{n\geq 0}\in C_b(\calh)$
is said to be $\calk$-convergent to a function $f\in C_b(\calh)$ (and we shall write
$f_n\overset{\calk}{\rightarrow} f$ or $f=\calk-\lim_{n\rightarrow\infty}f_n$) if for any compact set
$\calk\subset\calh$
$$
\sup_{n\in\N}\Vert f_n\Vert_\infty<+\infty \quad {\rm and } \quad
\lim_{n\rightarrow\infty}\sup_{x\in\calk}\vert f(x)-f_n(x)\vert =0.
$$
\noindent
Similarly, given $I\subseteq \R$, a sequence
$(f_n)_{n\geq 0}\in C_b(I \times \calh)$ is said to be $\calk$-convergent to a function $f\in C_b(I\times \calh)$ (and we shall write again
$f_n\overset{\calk}{\rightarrow} f$ or $f=\calk-\lim_{n\rightarrow\infty}f_n$) if for any compact set\footnote{In Section B.5.1 of \cite{FabbriGozziSwiech} the definition
of $\calk$-convergence on $I\times \calk$ is slightly different, as it requires
uniform convergence only on $I_0\times \calk$ for all compact subsets $I_0$ of $I$
and $\calk$ of $\calh$; it could be used here without affecting the results.}
$\calk\subset \calh$
$$
\sup_{n\in\N}\Vert f_n\Vert_\infty<+\infty \quad {\rm and } \quad
\lim_{n\rightarrow\infty}\sup_{(t,x)\in I\times\calk}\vert f(t,x)-f_n(t,x)\vert=0.
$$
\end{definition}


Now we recall the definition (given in  \cite{DP3}, beginning of Chapter 7)
of \textit{strict solution}
of a family of Kolmogorov equations.
Consider the following forward Kolmogorov equation with unknown $w$:
\begin{equation}\label{eq:Kforward}
  \left\{\begin{array}{l}\dis
\frac{\partial w(t,x)}{\partial t}=\frac{1}{2} Tr \;GG^*\; \nabla^2w(t,x)
+ \< Ax,\nabla w(t,x)\>
+\calf(t,x)
\qquad t\in [0,T],\,
x\in \calh,\\
\\
\dis w(0,x)=\phi(x).
\end{array}\right.
\end{equation}
where
the functions
$\calf:[0,T]\times \calh \to \R$ and $\phi:\calh\to\R$ are bounded and continuous.

\begin{definition}\label{df:strictandpistrong}
By \textit{strict solution} of the Kolmogorov equation (\ref{eq:Kforward}) we mean
a function $w$ such that
\begin{equation}\label{propr-strict-sol}
 \left\lbrace
\begin{array}{l}
w \in C_b([0,T]\times \calh)\quad  {\rm and}\quad w(0,x)=\phi(x)
\\
w(t,\cdot) \in UC^2_b(\calh), \; \forall t\in [0,T];
\\
w(\cdot, x)\in C^{1}([0,T]),\;\forall x\in D(A)\; \hbox{and $w$ satisfies (\ref{eq:Kforward}).}
\end{array}
\right.
\end{equation}
\end{definition}

Now we prove a key approximation lemma.

\begin{lemma}\label{lm:approximation}
Let Hypothesis \ref{ipotesibasic} and \ref{ipotesicostoconcreto} hold.
Let also (\ref{eq:hpdebreg}) or (\ref{eq:hpdebregbis}) hold.
Let $v$ be the mild solution of the HJB equation (\ref{HJBformale1})
and set $w(t,x)=v(T-t,x)$ for $(t,x)\in [0,T]\times \calh$.
Then there exist \red{two sequences of functions $(\bar\phi_n)$
and $(\calf_n)$ such that, for all $n\in \N$,
    \begin{equation}\label{eq:defapproxphif}
        \bar\phi_n \in C_c^\infty(\R^n), \qquad
        \calf_n  \in C_c^\infty([0,T]\times \calh) \cap \Sigma^1_{T,{1/2}}
    \end{equation}
and, calling $\bar\calf_n(t,x):=t^{1/2}\calf_n(t,x)$
    \begin{equation}\label{eq:convapproxphif}
        \bar\phi_n \rightarrow \bar\phi ,
        \qquad t^{1/2}\calf_n \to t^{1/2} H_{min} (\nabla^B w)
    \end{equation}}
in the sense of $\calk$-convergence, (the last in $(0,T]\times \calh$).
Moreover, defining \red{$\phi_n(x)=\bar\phi_n(x_0)$
and
\begin{equation}\label{v_n-mild}
w_n(t,x):= R_t \phi_n + \int_0^t\left[ R_{t-s}[\calf_n (s,\cdot) (x)+\ell_{0}(s)\right]ds
\end{equation}}
the following hold:
\begin{itemize}
  \item $w_n \in UC^{1,2}_b ([0,T]\times \calh) \cap \Sigma^1_{T,{1/2}}$,
  \item $w_n$ is a strict solution of
  \myref{eq:Kforward} with $\phi_n$ in place of $\phi$ and
\red{$\calf_n + \ell_{0}$ in place of $\calf $},
  \item we have, in the sense of $\calk$-convergence
  (the first in $[0,T]\times \calh$,
  the second in $(0,T]\times \calh$),
  \begin{equation}\label{eq:convapproxphifsol}
        w_n \rightarrow w ,
        \qquad t^{{1/2}}\nabla^B w_n \to t^{{1/2}}\nabla^B w.
    \end{equation}
\end{itemize}
\end{lemma}
\dim
We divide the proof in three steps.

{\bf Step 1: choosing the three approximating sequences.}
\red{We choose $\bar\phi_n$
to be the standard approximation by convolution of $\bar \phi$.}
To define $\calf_n$ we observe first that, since $w \in \Sigma^1_{T,{1/2}}$, then the function
$$
(t,x) \mapsto  \calf(t,x):= H_{min} (\nabla^B w(t,x))
$$
has the property that there exist $f:[0,T)\times \R^n \to \R$, continuous and bounded, such that
$$
\calf (t,x)=t^{-{1/2}} f(t,(e^{tA}x)_0).
$$
We then let $f_n$ be the approximation by convolution of $f$ and define
$$
\calf_n (t,x)=t^{-{1/2}} f_n(t,(e^{tA}x)_0).
$$
{\bf Step 2: proof that $w_n\in UC^{1,2}_b([0,T]\times\calh) \cap \Sigma^1_{T,{1/2}}$ and
that it is a strict solution.}
The fact that $w_n\in\Sigma^1_{T,{1/2}}$ follows immediately from \myref{v_n-mild},Proposition \ref{lemmaderhpdeb} and Lemma 5.4-(i) in \cite{FGFM}. Differentiability with respect to the variable $x$ follows by applying the dominated convergence theorem, while explicitely differentiating $ R_{t}\left[\phi\right]  $, namely or any $h\in\calh$
\begin{align*}
 \<\nabla R_{t}\left[\phi_n\right]  (x),h\>_H
&=\lim_{\alpha\rightarrow 0}\frac{1}{\alpha}
 \left[\int_\calh\phi_n\left(z+e^{tA}(x+\alpha h)\right)
 \caln\left(0,Q_{t}\right)  (dz)
 -\int_\calh
\phi_n\left(z+e^{tA}x\right)  \caln\left(  0,Q_{t}\right)  (dz) \right]
 \\[2mm]
 &  =\int_\calh\lim_{\alpha\rightarrow 0}\frac{1}{\alpha}
 \left[\phi_n\left(z+e^{tA}(x+\alpha h)\right)-\phi_n\left(z+e^{tA}x\right)\right]
 \caln\left(0,Q_{t}\right)  (dz)\\[2mm]
 &  =\int_\calh\<\nabla \phi_n\left(z+e^{tA}x\right), e^{tA}h\>_H
 \caln\left(0,Q_{t}\right)  (dz)=R_{t}\left[\<\nabla \phi_n , e^{tA}h\>_H\right ]  (x).
\end{align*}
In a similar way, differentiating twice
we get that, for all $h,k \in \calh$
\begin{align*}
& \<\nabla^2 R_{t}\left[\phi_n\right]  (x)h,k\>_\calh=
R_{t}\left[\<\nabla^2 \phi_n e^{tA}h,e^{tA}k\>_\calh\right ]  (x).
\end{align*}
Similarly we have, for the convolution term containing $\calf_n$,
\begin{align*}
& \<\nabla \int_0^t R_{t-s}\left[\calf_n(s,\cdot)\right]  (x),h\>_\calh
\\[2mm]
 &  =\lim_{\alpha\rightarrow 0}\frac{1}{\alpha}
 \int_0^t\left[\int_{\calh}\calf_n\left(s,z+e^{(t-s)A}(x+\alpha h)\right)
 \caln\left(0,Q_{t-s}\right)  (dz)
 -\int_{\calh}
 \calf_n\left(s,z+e^{(t-s)A}x\right)  \caln\left(  0,Q_{t-s}\right)  (dz) \,ds\right]
 \\[2mm]
 &  =\int_0^t\int_{\calh}\lim_{\alpha\rightarrow 0}\frac{1}{\alpha}
 \left[\calf_n\left(s,z+e^{(t-s)A}(x+\alpha h)\right)-\calf_n\left(s,z+e^{(t-s)A}x\right)\right]
 \caln\left(0,Q_{t-s}\right)  (dz)\,ds\\[2mm]
 &  =\int_0^t\int_{\calh}\<\nabla \calf_n\left(s,z+e^{(t-s)A}x\right), e^{(t-s)A}h\>_\calh
 \caln\left(0,Q_{t-s}\right)  (dz)\,ds\\[2mm]
 &=\int_0^tR_{t-s}\left[\<\nabla \calf_n(s,\cdot) ,e^{(t-s)A}h\>_\calh\right ]  (x)\,ds,
\end{align*}
and also, arguing in the same way,
\begin{align*}
& \<\nabla^2\int_0^t R_{t}\left[\calf_n(s,\cdot)\right]  (x)\,ds\, h,k\>_\calh=
\int_0^tR_{t-s}\left[\<\nabla^2 \calf_n(s,\cdot) e^{(t-s)A}h,e^{(t-s)A}k\>_\calh\right ]  (x)\,ds.
\end{align*}
The  convolution term involving $\ell_{0,n}$ is treated exactly in the same way.
The proof that $w_n$ is differentiable with respect to time and that
$w_{nt}\in UC_b([0,T]\times\calh) $ is completely analogous
to what is done in \cite[Theorems 9.23 and 9.25]{DP1})
for homogeneous Kolmogorov equations and we omit it\footnote{The proof
in the nonhomogeneous case can be found in the recent book \cite[Theorem 4.135, Step 2, or Proposition B.91]{FabbriGozziSwiech}}.
By Theorem 5.3 in \cite{CeGo}, see also Theorem 7.5.1 in \cite{DP3} for Kolmogorov equations,
we finally conclude that $w_n$
is a strict solution to equation \ref{eq:Kforward}.

{\bf Step 3: proof of the convergences.}
First we prove that the sequences $(w_n)$ and $(t^{1/2}\nabla^B w_n)$
are bounded uniformly with respect to $n$. Indeed,
by (\ref{v_n-mild}) and by the properties of convolutions,
\begin{align*}
&\vert w_n(t,x)\vert
\leq \Vert \bar\phi_n\Vert_\infty +
\int_0^t\sup_{x\in\calh}\left[\vert\calf_n (s,x)\vert+\vert
\red{\ell_{0}(s)}\vert\right]\,ds
\\
&\leq
\Vert \bar\phi\Vert_\infty + \int_0^t\sup_{y\in\R^n}\left[s^{-{1/2}} \vert
f_n (s,y)\vert+\vert
\red{\bar\ell_{0}(s)}\vert\right]\,ds \leq
\Vert \bar\phi\Vert_\infty + \int_0^t\sup_{y\in\R^n}\left[ s^{-{1/2}}\vert f (s,y)\vert+\vert\red{\bar\ell_{0}(s)}\vert\right]\,ds
\end{align*}
Moreover, using Proposition \ref{lemmaderhpdeb} and the results in \cite{FGFM}, Section 5, formula 5.9,  with $\psi=identity$, namely
\begin{align}
&\<\nabla ^B \left(\int_{0}^{t}
R_{t-s}\left[\nabla^{B}(g(s,\cdot)\right]  ds
\right) (x),k \>_{\R^m}=
\label{eq:derBconvnew}
\\[3mm]
\nonumber
&=
\int_{0}^{t} \int_\calh
s^{-{1/2}} \bar f\left(s, (e^{sA}z)_0+(e^{tA}x)_0\right)
\<(Q^0_{t-s})^{-{1/2}}
\left( e^{tA} Bk\right)_0, (Q^0_{t-s})^{-{1/2}}z_0\>_{\R^n}
\caln\left(0,Q_{t-s}\right)(dz)ds,
\end{align}
we get
\begin{equation*}
 \vert t^{1/2} \nabla^B w_n(t,x)\vert
 \leq C\Vert \bar\phi\Vert_\infty
 +Ct^{1/2}\int_0^ts^{-{1/2}}(t-s)^{-{1/2}}\sup_{y\in\R^n}\left(\vert
   f\left(s, y\right)\vert+  \vert\red{\bar\ell_{0}(s)}\vert\right)  ds,
\end{equation*}
for a suitable $C>0$.
Now, with similar computations, we prove the convergences. Indeed,
\begin{align*}
 \vert w_n(t,x)&-w(t,x)
 \vert\leq \int_\calh\vert \bar\phi_n\left( (e^{tA}z)_0+(e^{tA}x)_0\right)-
\bar\phi\left( (e^{tA}z)_0+(e^{tA}x)_0\right)\vert\caln(0,Q_{t})(dz) \\
&\red{ +
\int_0^t\int_\calh\left[\vert s^{-{1/2}} \calf_n \left(s,(e^{sA}z)_0+(e^{tA}x)_0\right)
-
s^{-{1/2}} \calf \left(s,(e^{sA}z)_0+(e^{tA}x)_0\right)\vert
\right]\caln(0,Q_{t-s})(dz)\,ds.}
\end{align*}
Since, for every compact set $\calk \subset \calh$ the set
$\{(e^{tA}x)_0, \; t \in [0,T],\, x \in \calk\}$ is compact in $\R^n$,
then by the Dominated Convergence Theorem we get that for any compact set $\calk \subset \calh$
\begin{equation}
\label{eq:convvnnew}
\sup_{(t,x)\in[0,T]\times\calk} \vert w_n(t,x)-w(t,x)\vert\rightarrow 0.
\end{equation}
Moreover, using again Proposition \ref{lemmaderhpdeb}
and Lemma 5.4 in \cite{FGFM},
for a suitable $C>0$, we get\red{
\begin{align*}
 &\vert  \nabla^B w_n(t,x)- \nabla^B w(t,x)\vert\\
&\leq C
\int_\calh\left( \bar\phi_n\left( (e^{tA}z)_0+(e^{tA}x)_0\right)-
\bar\phi\left( (e^{tA}z)_0+(e^{tA}x)_0\right)\right)\<(Q^0_{t})^{-{1/2}}
\left( e^{tA} Bh\right)_0, (Q^0_{t})^{-{1/2}}z_0\>_{\R^n}\caln(0,Q_{t})(dz)\\
  & +C\int_0^t
  \int_\calh  \left(\left\vert
 s^{-{1/2}}  f_n\left(s, (e^{sA}z)_0+(e^{tA}x)_0\right)
 -s^{-{1/2}}  f\left(s, (e^{sA}z)_0+(e^{tA}x)_0\right)
 \right\vert\right)\\
 &\quad \left \vert
 \<(Q^0_{t-s})^{-{1/2}}
 \left( e^{tA} Bh\right)_0, (Q^0_{t-s})^{-{1/2}}z_0\>_{\R^n}
  \right\vert\caln(0,Q_{t-s})(dz) ds.
\end{align*}}
By Proposition \ref{lemmaderhpdeb} we know that for suitable $C>0$,
\[
 \int_\calh\vert\<(Q^0_{t})^{-{1/2}}
\left( e^{tA} Bh\right)_0, (Q^0_{t})^{-{1/2}}z_0\>_{\R^n}\vert^2\caln(0,Q_{t})(dz)
\leq C\Vert (Q^0_{t})^{-{1/2}}
\left( e^{tA} Bh\right)_0\Vert^2\leq \frac{C}{t}.
\]
Hence, applying Cauchy-Schwartz inequality we get\red{
\begin{align*}
&\vert  \nabla^B w_n(t,x)- \nabla^B w(t,x)\vert\\
&\leq C t^{-{1/2}}
 \vert\int_\calh\vert \bar\phi_n\left( (e^{tA}z)_0+(e^{tA}x)_0\right)-
 \bar\phi_n\left( (e^{tA}z)_0+(e^{tA}x)_0\right)\vert^2\caln(0,Q_{t})(dz)\vert^{1/2}\\
  &+\int_0^t (t- s)^{-{1/2}}
 \left( \int_\calh \left\vert
 s^{-{1/2}}  f_n\left(s, (e^{sA}z)_0+(e^{tA}x)_0\right)
 -s^{-{1/2}}  f\left(s, (e^{sA}z)_0+(e^{tA}x)_0\right)
 \right\vert
 \caln(0,Q_{t-s})(dz)\right)^{1/2}.
\end{align*}}
Applying the Dominated Convergence Theorem as for the proof of
\myref{eq:convvnnew} we get the final claim
\[
\sup_{x\in(0,T]\times\calk} \vert t^{1/2}\nabla^B w_n(t,x)-t^{1/2}\nabla^B w(t,x)\vert\rightarrow 0, \qquad  \hbox{for any compact set $\calk\subset\calh$.}
\]
\qed
\medskip

Notice that, using the terminology of \cite{CeGo,G1}, the above result
implies that a mild solution (\ref{eq:Kforward}) is also a {\em $\calk$-strong solution}.
In general, in an infinite dimensional Hilbert space $H$, existence of $\calk$-strong solutions is not a routine application of the theory of evolution equations, as the operator $\call$ formally introduced in \ref{eq:ell}
is not the infinitesimal generator of a strongly continuous semigroup in the Banach space $C_b(H)$. To overcome this difficulty in the already mentioned paper
\cite{CeGo} the theory of weakly continuous (or $\calk$-continuous)
semigroups has been used.

\begin{remark}\label{rm:defpistrictstrong}
The approximation results proved just above is needed to prove the fundamental identity,
(see next Proposition \ref{prop rel fond}) which is the key point to get the verification theorem and the existence of optimal feedback controls.
The idea is to apply Ito's formula to the approximating sequence $w_n$ composed with
the state process $Y$ and then to pass to the limit for $n\to +\infty$
(see e.g. \cite{G1} or \cite[Section 4.4]{FabbriGozziSwiech}).
However in the literature the approximating sequence is taken more regular,
i.e. the $w_n$ are required to be classical solutions
(see e.g. \cite[Section 6.2, p.103]{DP3}) of \myref{eq:Kforward}.
This in particular means that $\nabla w_n\in D(A^*)$ and this fact is crucial
since it makes possible to apply Ito's formula without requiring that the
state process $Y$ belongs to $D(A)$, which would be a too strong requirement, see fore related results \cite{FlaZan}.

In our case the used approximating procedure does not give rise in general to
functions $w_n$ with $\nabla w_n\in D(A^*)$.
Indeed for our purposes we need that the approximants of the data $\phi,\,\ell_0, \,\calf$ remain all in the space $\Sigma^1_{T,{1/2}}$; without this, since we only have
``partial'' smoothing,
it is not clear at all how to prove the convergence of the derivative $\nabla^B w_n$
(which is needed when we pass to the limit to prove the fundamental identity in next subsection).
Hence, in particular, since we need that, for all $n\in\N$, $\calf_n\in\Sigma^1_{T,{1/2}}$, and since
$\calf$ is written in terms of $f$, we approximate $f$
by $f_n$, and this procedure gives the approximants $\calf_n$ of $\calf$.
In this way $\calf_n\in \Sigma^1_{T,{1/2}}$ but $\nabla \calf_n \notin D(A^*)$.

Summing up, we are only able to find approximating strict solutions
and not classical solutions. Since the state process $Y$ does not belong to $D(A)$
this fact will force us to introduce suitable regularizations $Y_k$ of it
(see the proof of Proposition \ref{prop rel fond}).
\end{remark}

\begin{remark}\label{rm:defpistrictstrongbis}
Calling \red{$\ell_n:=\ell_{0}+\calf_n -H_{min}(\nabla^Bw_n)$}
is not difficult to see that the sequence
$w_n$ is a strict solution of the approximating HJB equation
\begin{equation*}
  \left\{\begin{array}{l}\dis
\frac{\partial w(t,x)}{\partial t}=
\frac{1}{2}Tr \;GG^*\nabla^2w(t,x)
+ \< Ax,\nabla w(t,x)\>_\calh
+H_{min} (\nabla^B w(t,x)) +\ell_n(t,x),\qquad t\in [0,T],\,
x\in \calh,\\
\\
\dis w(0,x)=\phi_n(x),
\end{array}\right.
\end{equation*}
This means, with the terminology used e.g. in \cite{G1}, that $w$ is a $\calk$-strong solution
of \myref{HJBINTRO}. We do not go deeper into this since here we use the approximation
only as a tool to solve our stochastic optimal control problem.
\end{remark}

\begin{remark}\label{rm:piconv}
It can be noticed, see also \cite{FGFM}, that our results on the HJB equation
could be extended without difficulties to the case
when the boundedness assumption on $\bar \phi$ and $\bar \ell_0$, and so on  $\phi$ and $\ell_0$,
is replaced by a polynomial growth assumption:
namely that, for some $N\in \N$, \red{the function
\begin{equation}\label{eq:polgrowthphil0}
 x\mapsto \dfrac{\phi(x)}{1+\vert x\vert ^N},
\end{equation}
is} bounded.
Also \ref{lm:approximation} can be easily generalized to the case when the data $\phi$
and $\ell_0$ are not bounded but satisfy a polynomial growth condition in the variable $x$
as from \myref{eq:polgrowthphil0}.

Moreover it is still possible to extend the results of Lemma \ref{lm:approximation} to the case when $\phi$ \red{ is} only measurable. In this case the approximations would take place in the sense of the $\pi$-convergence, which is weaker than the $\calk$-convergence and towards the $\calk$-convergence has also the disadvantage of being not metrizable. For more on the notion of $\pi$-convergence the reader can see
\cite[Section 6.3]{DP3} (also \cite{EthierKurtz}, \cite[Section B.5]{FabbriGozziSwiech} and \cite{PriolaStudia}).
\end{remark}

\section{Verification Theorem and Optimal Feedbacks}\label{sec-verifica}

The aim of this section is to provide a verification theorem and the existence of optimal feedback controls for our problem.
This in particular will imply that the mild solution $v$ of the HJB equation (\ref{HJBINTRO}) built in Theorem \ref{esistenzaHJB} is equal to the
value function $V$ of our optimal control problem.

The main tool needed to get the wanted results is an identity (often called ``{\em fundamental identity}'', see equation (\ref{relfond})) satisfied by the solutions of the HJB equation. When the solution is smooth enough
(e.g. it belongs to $UC_{b}^{1,2}\left(  \left[0,T\right]  \times H\right)$)
such identity is easily proved using the Ito's Formula. Since in our case the value
function does not possess this regularity, we proceed by approximation, following the lines of Section \ref{sec-Kstrong}.
Due to the features of our problem (lack of smoothing and of the structure condition)
the methods of proof used in the literature do not apply here.
We will discuss the main issues along the way.

\subsection{The Fundamental Identity and the Verification Theorem}

Now we finally go back to the control problem of Section \ref{section-statement}.
We rewrite here for the reader convenience the state equation (\ref{eq-contr-rit}),
\begin{equation*}
\left\{
\begin{array}
[c]{l}
dy(s)  =a_0 y(s) ds+b_0 u(s) ds +\int_{-d}^0b_1(\xi)u(s+\xi)d\xi+\sigma dW_s
,\text{ \ \ \ }s\in[t,T] \\
y(t)  =y_0,\\
u(\xi)=u_0(\xi), \quad \xi \in [-d,0),
\end{array}
\right.
\end{equation*}
and its abstract reformulation (\ref{eq-astr}),
\begin{equation*}
\left\{
\begin{array}
[c]{l}
dY(s)  =AY(s) ds+Bu(s) ds+GdW_s
,\text{ \ \ \ }s\in[t,T] \\
Y(t)=x=(x_0,x_1).
\end{array}
\right.
\end{equation*}
Similarly the cost functional in (\ref{costoconcreto}) is
\begin{equation*}
J(t,x;u)=\E \int_t^T \left(\red{\bar\ell_0(s)}+\ell_1(u(s))\right)\;ds +\E  \bar\phi(x(T))
\end{equation*}
and is rewritten as (see (\ref{costoconcreto1}))
\begin{equation*}
J(t,x;u)=\E \int_t^T \left(\red{\ell_0(s)}+\ell_1(u(s))\right)\;ds +\E  \phi(Y(T)).
\end{equation*}
We notice that throughout this subsection and the following one, in order to avoid further technical difficulties,  we keep the probability space $(\Omega,\calf,\P)$
fixed. Nothing would change if we work in the weak formulation, where the probability space can change (see e.g. \cite{YongZhou99}[Chapter 2] and
\cite{FabbriGozziSwiech}[Chapter 2] for more on strong and weak formulations
in finite and infinite dimension, respectively).
We first prove the fundamental identity.
\begin{proposition}\label{prop rel fond}
Let Hypotheses \ref{ipotesibasic} and \ref{ipotesicostoconcreto} hold.
Let also (\ref{eq:hpdebreg}) or (\ref{eq:hpdebregbis}) hold.
Let $v$ be the mild solution of the HJB equation (\ref{HJBINTRO})
according to Definition \ref{defsolmildHJB}.
Then for every $t\in[ 0,T] $ and $x\in
H$, and for every admissible control $u$, we
have the fundamental identity
\begin{equation}\label{relfond}
 v(t,x)
=J(t,x;u)+\E\int_t^T \left[H_{min}(\nabla^B v(s,Y(s)))
- H_{CV}(\nabla^B v(s,Y(s));u(s))\right]\,ds.
\end{equation}
\end{proposition}
\dim
Take any admissible state-control couple $(Y(\cdot),u(\cdot))$,
and let $v_n(t,x):=w_n(T-t,x)$ where $(w_n)_n$ is the approximating sequence
of strict solutions defined
in Lemma \ref{lm:approximation}. We want to apply the Ito formula to
$v_n(t, Y(t))$. Unfortunately, as mentioned in Remark \ref{rm:defpistrictstrong},
this is not possible since the process
$Y(t)$ does not live in $D(A)$. So we approximate it as follows.
Set, for $k \in \N$, sufficiently large,
\begin{equation}\label{eq:Ykdef}
Y_k(s;t,x)=k(k-A)^{-1}Y(s;t,x).
\end{equation}
The process $Y_k$ is in $D(A)$, it converges to $Y$ ($\P$-a.s. and $s\in [t,T]$ a.e.)
and it is a strong solution\footnote{Here we mean strong in the probabilistic sense
and also in the sense of \cite{DP1}, Section 5.6.}
of the Cauchy problem
\begin{equation*}
\left\{
\begin{array}
[c]{l}
dY_k(s)  =AY_k(s) ds+B_ku(s) ds+G_kdW_s
,\text{ \ \ \ }s\in [t,T] \\
Y_k(t)=x_k,
\end{array}
\right.
\end{equation*}
where $B_k=k(k-A)^{-1}B$, $G_k=k(k-A)^{-1}G$ and $x_k=k(k-A)^{-1}x$.
Now observe that, thanks to \myref{resolvent} and \myref{B},
the operator $B_k$ is continuous, hence we can apply Dynkin's formula (see e.g.
\cite[Section 1.7]{FabbriGozziSwiech} or \cite[Section 4.5]{DP1})
to $v_{n}(s, Y_k(s))$ in the interval $[t,T]$, getting
\begin{align}\label{Dynkinv^nk}
&\E v_n(Y_k(T)) - v_n(t,x_k)\\ \nonumber &=
\E\int_t^T
\left[v_{nt}(s,Y_k(s))+  \frac{1}{2}Tr \;GG^*\nabla^2v(s,Y_k(s))
+ \< AY_k(s)+B_ku(s),\nabla v_n(s,Y_k(s))\>_\calh\right]ds.
\end{align}
Using the Kolmogorov equation (\ref{eq:Kforward}), whose strict solution is $w_n$,
we then write
\begin{equation}\label{quasirelfondv^nk}
\E\phi_n(Y_k(T)) - v_n(t,x_k)=\E\int_t^T \left[\calf_n (s,Y_k(s))+ \ell_{0,n}(s,Y_k(s))+\<B_k u(s),\nabla v_n(s,Y_k(s))\>_{\calh} \right]ds
\end{equation}
We first let $k\rightarrow\infty$ in \myref{quasirelfondv^nk}.
Since $\ell_{0,n}$ and $\nabla v_n$ are
bounded functions and since $\calf_n(s,x)$ has a singularity of type $s^{-{1/2}}$ with
respect to time and is bounded with respect to $x$, we can apply
the Dominated Convergence Theorem to all terms but the last getting
\begin{align}\label{quasirelfondv^n}
&\E\phi_n(Y(T)) - v_n(t,x)\\ \nonumber&=\E\int_t^T [\calf_n (s,Y(s))+\red{\ell_{0}(s)}]ds
+\lim_{k \to + \infty}\E\int_t^T\<B_k u(s),\nabla v_n(s,Y_k(s))\>_\calh]ds.
\end{align}
Concerning the last term we observe first that, by \myref{resolvent} and
 \myref{B},  as well as by property of the operators $k(k-A)^{-1}$, we have that for every $u\in\R^m$ $$B_ku \to Bu.$$
This in particular implies, by the Banach-Steinhaus Theorem, that
$\{B_k u\}_k$ is uniformly bounded in $\calh$.
Now we use the fact that $\nabla v_n(s,x) \in \calh$
(see \myref{eq:nablaperSigma})
to rewrite the integrand of the last term of \myref{quasirelfondv^n} as
$$
\<B_k u(s),\nabla v_n(s,Y_k(s))-\nabla v_n(s,Y(s))\>_{\calh}
+\<B_k u(s),\nabla v_n(s,Y(s))\>_{\calh}
$$
Thanks to what said above the first term goes to $0$ as $k \to + \infty$ and is dominated while the second term is also dominated and converges to
$\<B u(s),\nabla v_n(s,Y(s))\>_{\calh}$
which, thanks to \myref{eq:nablaperSigma}, is equal to
$\< u(s),\nabla^B v_n(s,Y(s))\>_{\R^m}$ (both convergences are clearly
$\P$-a.s. and $s\in [t,T]$ a.e.). Hence
$$
\lim_{k \to + \infty}\E\int_t^T\<B_k u(s),\nabla v_n(s,Y_k(s))\>_\calh ds
=
\E\int_t^T\<u(s),\nabla^B v_n(s,Y(s))\>_{\R^m}ds
$$
Now we let $n\rightarrow\infty$. By Lemma \ref{lm:approximation}, we know that
$$v_n(t,x)\rightarrow v(t,x) \quad\text{  and  }\quad(T-t)^{1/2} \nabla^Bv_n(t,x)\rightarrow (T-t)^{1/2}\nabla^B v(t,x)$$
pointwise. Moreover $v_n(t,x)$, $(T-t)^{1/2}\nabla^Bv^n(t,x)$
are uniformly bounded, so that, by dominated convergence, we get
\[
\E\int_t^T \<u(s),\nabla^Bv_n(s,Y(s))\>_{\R^m}]\,ds\rightarrow
\E\int_t^T \<u(s),\nabla^Bv(s,Y(s))\>_{\R^m}]\,ds.
\]
The convergence
\[
 \E\phi_n(x(T)) -\E\int_t^T  [\calf_n(s,Y(s))+ \red{\ell_{0}(s)}]ds
 \rightarrow
\E\phi(x(T)) -\E\int_t^T [H_{min}(\nabla^B v(s,Y(s)))+\red{ \ell_0(s)}]ds
\]
follows directly by the construction of the approximating sequences \red{$(\phi_n)_n$ and $(\calf_n)_n$}.
Then, adding and subtracting $\E\dis\int_t^T  \ell_1(u(s))ds$
and letting $n\rightarrow\infty$ in (\ref{quasirelfondv^n})
we obtain
\begin{align*}
 v(t,x)&
=\E\phi(Y(T))+\E\int_t^T [\red{\ell_0(s)}+ \ell_1(u(s))]ds\\
&+\E\int_t^T \left[H_{min}(\nabla^B v(s,Y(s)))-
H_{CV}(\nabla^B v(s,Y(s));u(s))\right]
\,ds
\end{align*}
which immediately gives the claim.
\qed


\begin{remark}
\label{rm:Yk}
One may wonder why we approximate the process $Y$ with $Y_k$ as in \myref{eq:Ykdef}
instead of using the Yosida approximants $A_k$ of $A$ as, e.g., in
\cite[p.144]{DP3}. The reason is that we need that $Y_k$ belongs to $D(A)$,
which is not guaranteed if we use Yosida approximants.
A similar procedure is used, in a different context, in
the book \cite{DP1}, in the proof of Theorem 7.7, p. 203.
\end{remark}

We can now pass to prove our Verification Theorem i.e.
a sufficient condition of optimality given in term of the
mild solution $v$ of the HJB equation (\ref{HJBINTRO}).

\begin{theorem}
\label{teorema controllo}
Let Hypotheses \ref{ipotesibasic} and \ref{ipotesicostoconcreto} hold.
Let also (\ref{eq:hpdebreg}) or (\ref{eq:hpdebregbis}) hold.
Let $v$ be the mild solution of the HJB equation (\ref{HJBINTRO})
whose existence and uniqueness is proved in (\ref{esistenzaHJB}).
Then the following holds.
\begin{itemize}
 \item For all $(t,x)\in [0,T]\times \calh$ we have
$v(t,x) \le V(t,x)$, where $V$ is the value function
defined in (\ref{valuefunction}).
\item
Let $t\in [0,T]$ and $x\in \calh$ be fixed.
If, for an admissible control $u^*$, we
have, calling $Y^*$ the corresponding state,
$$
u^*(s)\in \arg\min_{u\in U}H_{CV}(\nabla^B v(s,Y^*(s);u)
$$
for a.e. $s\in [t,T]$, $\P$-a.s., then the pair $(u^*,Y^*)$ is
optimal for the control problem starting from $x$ at time $t$
and $v(t,x)=V(t,x)=J(t,x;u^*)$.
\end{itemize}
\end{theorem}
\dim The first statement follows directly by \eqref{relfond} due to the
negativity of the integrand.
Concerning the second statement, we immediately see that, when $u=u^*$
(\ref{relfond}) becomes
$v(t,x)=J(t,x;u^*)$.
Since we know that for any admissible control $u$
\[
 J(t,x;u)\geq V(t,x) \ge v(t,x),
\]
the claim immediately follows.
\qed

\subsection{Optimal feedback controls and $v=V$}
\label{sec:contr-feedback}
We now prove the existence of optimal feedback controls. Under the Hypotheses
of Theorem \ref{teorema controllo} we define, for $(s,y)\in [0,T)\times \calh$,
the {\em feedback map}
\begin{equation}\label{defdiPsi}
\Psi(s,y):=\arg \min_{u\in U} H_{CV}(\nabla^Bv(s,y);u),
\end{equation}
where, as usual, $v$ is the solution of the HJB equation \myref{HJBINTRO}.
Given any $(t,x)\in [0,T)\times \calh$,
the so-called Closed Loop Equation (which here is, in general, an inclusion)
is written, formally, as
\begin{equation}\label{cleinclusion}
\left\{
\begin{array}
[c]{l}
dY(s) \in AY(s) ds+B\Psi\left(s,Y(s)\right) ds+GdW_s
,\text{ \ \ \ }s\in [t,T) \\
Y(t)=x.
\end{array}
\right.
\end{equation}
First of all we have the following straightforward corollary whose proof is immediate from Theorem \ref{teorema controllo}.

\begin{corollary}
\label{cr:optimalfeedback}
Let the assumptions of Theorem \ref{teorema controllo} hold true.
Let $v$ be the mild solution of \myref{HJBformale1}.
Fix $(t,x)\in [0,T)\times \calh$ and assume that, on $[t,T)\times \calh$, the map $\Psi$
defined in (\ref{defdiPsi}) admits a measurable selection
$\psi:[t,T)\times \calh\to \Lambda$ such that the Closed Loop Equation
\begin{equation}
\label{eq:CLEselection}
\left \{
\begin{array}{l}
d Y(s) = AY(s)d s+B\psi\left(s,Y(s)\right) ds+GdW_s
,\text{ \ \ \ }s\in [t,T) \\
Y(t)=x.
\end{array}
\right.
\end{equation}
has a mild solution $Y_\psi(\cdot;t,x)$  (in the sense of \cite[p.187]{DP1}).
Define, for $s \in [t,T)$, $u_\psi (s)=\psi(s,Y_\psi(s;t,x))$.
Then the couple
$(u_\psi(\cdot),Y_\psi(\cdot;t,x))$ is optimal at
$(t,x)$ and $v(t,x)=V(t,x)$.
If, finally, $\Psi(t,x)$ is always a singleton and the mild solution
of \myref{eq:CLEselection} is unique,
then the optimal control is unique.
\end{corollary}

We now give sufficient conditions to verify the assumptions of Corollary
\ref{cr:optimalfeedback}.
First of all define
 \begin{equation}\label{defdigammagrandebis}
\Gamma(p):=\left\{ u\in U: \<p,u\>+\ell_1(u)= H_{min }(p)\right\}.
\end{equation}
Then, clearly, we have $\Psi(t,x)=\Gamma(\nabla^B v(t,x))$.
Observe that, under mild additional conditions on $U$ and $\ell_1$
(for example taking $U$ compact or $\ell_1$ of superlinear growth),
the set $\Gamma$ is nonempty for all $p \in \R^m$.
If this is the case then, by \cite{AubFr}, Theorems 8.2.10 and
8.2.11, $\Gamma$ admits a measurable selection, i.e. there exists
a measurable function $\gamma:\R^m \rightarrow U$ with
$\gamma(z)\in \Gamma(z)$ for every $z\in \R^m$.
Since $H_{min }$ is Lipschitz continuous, then $\Gamma$, and so
$\gamma$, must be uniformly bounded. In some cases
studied in the literature this is enough to find an optimal feedback
but not in our case (read on this the subsequent Remark \ref{rm:discfeed}-(ii)).
Hence to prove existence of a mild solution of the closed loop
equation \myref{eq:CLEselection}, as requested in Corollary
\ref{cr:optimalfeedback}, we need more regularity of the feedback
term $\psi(s,y)=\gamma(\nabla^B v(s,y))$.
Beyond the smooth assumptions on the coefficients required in the second part of Theorem \ref{esistenzaHJB}, which give the regularity of
$\nabla^B v(t,x)$, we need the following assumption about the map $\Gamma$.

\begin{hypothesis}\label{hp:lipsel}
The set-valued map $\Gamma$ defined in (\ref{defdigammagrandebis})
admits a Lipschitz continuous selection $\gamma$.
\end{hypothesis}

\begin{remark}

\label{rm:discfeed}
The problem of the lack of regularity of the feedback law is sometimes faced (see e.g. in \cite{FT2}) by formulating the optimal control problem in the weak sense (see e.g. \cite{FlSo}, or \cite[Proposition 4.199 and Theorem 6.36]{FabbriGozziSwiech}, or \cite{YongZhou99}, Section 4.2)
and then using Girsanov Theorem to prove existence, in the weak sense, of a mild solution of (\ref{eq:CLEselection}) when the map $\psi$ is only measurable and bounded.
This is not possible here due to the already mentioned absence of the structure condition in the controlled state equation (i.e. the control affects the system not only through the noise).
\end{remark}

Taking the selection $\gamma$ from Hypothesis \ref{hp:lipsel}
we consider the closed loop equation
\begin{equation}\label{cle}
\left\{
\begin{array}
[c]{l}
dY(s)  =AY(s) ds+B\gamma\left(\nabla^B v(s,Y(s))\right) ds+GdW_s
,\text{ \ \ \ }s\in[t,T] \\
Y(t)=x=(x_0,x_1),
\end{array}
\right.
\end{equation}
and we have the following result.
\begin{theorem}\label{teo su controllo feedback}
Assume that Hypotheses \ref{ipotesibasic}, \ref{ipotesicostoconcreto}, \ref{ipotesicostoconcretobis} and \ref{hp:lipsel} hold true.
Fix any $(t,x)\in [0,T)\times H$.
Let also (\ref{eq:hpdebreg}) or (\ref{eq:hpdebregbis}) hold and let either $\bar\phi$ be Lipschitz continuous or $\sigma$ be invertible.
Then the closed loop equation (\ref{cle}) admits a unique mild solution $Y_\gamma(\cdot;t,x)$ (in the sense of \cite[p.187]{DP1}) and setting
$$
u_\gamma(s)=\gamma\left(\nabla^{B} v(s,Y_\gamma(s;t,x) )\right), \quad s \in [t,T]
$$
we obtain an optimal control at $(t,x)$ which is unique if $\Gamma$ is always a singleton. Moreover $v(t,x)=V(t,x)$.
\end{theorem}

\noindent {\bf Proof.} Thanks to Corollary \ref{cr:optimalfeedback}
it is enough to prove the existence and uniqueness of the mild solution of (\ref{cle}). We apply a fixed point theorem to the following integral form of
(\ref{cle}):
\begin{equation}\label{clemild}
Y(s)= e^{(s-t)A}x+
\int_t^s e^{(s-r)A}G\,dW_r
+\int_{t}^s e^{(s-r)A}B\gamma(\nabla^{B}v(r ,\overline{X}_r)))dr.
\end{equation}
By Hypothesis \ref{ipotesicostoconcretobis} and the second part of Theorem \ref{esistenzaHJB} we get that, when $\sigma$ is invertible, the mild solution $v$ of the HJB equation (\ref{HJBformale1}) is differentiable, with bounded derivative. Moreover, since, again by the second part of
Theorem \ref{esistenzaHJB},
$v$ admits the second order derivative $\nabla^2 v$, 
we deduce that $t^{1/2}\nabla^B v(t,\cdot)$ is
Lipschitz continuous, uniformly with respect to $t$.
Such Lipschitz property can be proved also when $\bar\phi$ is Lipschitz continuous by taking $\nabla^B$ of both sides of (\ref{solmildHJB}) and applying Gronwall inequality.
Using this Lipschitz property we can solve (\ref{clemild}) by a fixed point argument. Since the argument to do this is straightforward
we only show how to estimate the difficult term in \myref{clemild}. We have
\begin{align*}
\int_{t}^s& \vert e^{(s-r)A}B\left(\gamma(\nabla^{B}
v(r ,\overline{X}(r)))-\gamma(\nabla^{B}
v(r ,\overline{Y}(r)))\right)\vert_\calh dr
\leq \int_{t}^s C\vert\gamma(\nabla^{B}
v(r ,\overline{X}(r)))-\gamma(\nabla^{B}
v(r ,\overline{Y}(r)))\vert_{\R^m}] dr
\\
&
\leq C\int_{t}^s \vert\nabla^{B}
v(r ,\overline{X}(r)))-\nabla^{B}
v(r ,\overline{Y}(r))\vert_{\R^m}] dr
\leq C\int_{t}^s r^{-{1/2}}\vert \overline{X}(r)-\overline{Y}(r)\vert_\calh dr
\end{align*}
where $C$ is a constant that can change its values from line to line.
\qed

\begin{corollary} Under the assumptions of Theorem \ref{teo su controllo feedback},
if we assume that $\ell_1:\R^m\rightarrow \R$ is differentiable with an invertible derivative and that $(\ell_1')^{-1}$ is Lipschitz continuous,
then the closed loop equation (\ref{cle}) admits a unique mild solution
$Y_\gamma(\cdot;t,x)$ (in the sense of \cite[p.187]{DP1}).
In such case, setting
$$
u_\gamma(s)=\gamma\left(\nabla^{B} v(s,Y_\gamma(s;t,x) )\right), \quad s \in [t,T]
$$
we get that $u_\gamma$ is the unique optimal control at $(t,x)$. Moreover $v(t,x)=V(t,x)$.
\end{corollary}
\dim
The proof is a straightforward consequence of Theorem \ref{teo su controllo feedback}. Indeed in this case, see (\ref{psi1}),
the infimum of $H_{CV}$ is achieved at
\[
u=(\ell_1')^{-1}(z), \text{ so that }\Gamma(z)=\{(\ell_1')^{-1}(z)\}.
\]
Since $\Gamma(z)$ is always a singleton we get uniqueness from Corollary \ref{cr:optimalfeedback}.
\qed

We devote our final result to show that the identification $v=V$ can be done,
using an approximation procedure, also
in cases when we do not know if optimal feedback controls exist.


\begin{theorem}\label{teo:v=V}
Let Hypotheses \ref{ipotesibasic}, \ref{ipotesicostoconcreto} hold.
Let also (\ref{eq:hpdebreg}) or (\ref{eq:hpdebregbis}) hold.
Moreover let Hypotheses \ref{ipotesicostoconcretobis}-(ii)
and \ref{hp:lipsel} hold and let
$\phi$ and $\ell_0$ be uniformly continuous.
Then $v=V$.
%
\end{theorem}
\dim
We approximate \red{
$\phi$ by approximating 
$\bar\phi$ with standard approximants 
$\bar\phi_n$ built by convolution}.
We set
\begin{equation}\label{costo_n}
J_n(t,x;u)=\E \int_t^T \left(\red{\ell_{0}(s)}+\ell_1(u(s))\right)\;ds +\E  \phi_n(Y(T))
\end{equation}
and call $w_n$ the mild solution of the HJB equation
\begin{equation}\label{HJBformale-n}
  \left\{\begin{array}{l}\dis
\frac{\partial w(t,x)}{\partial t}=\call [w(t,\cdot)](x) +\red{\ell_{0}(t)}+
H_{min} (\nabla^B w(t,x)),\qquad t\in [0,T],\,
x\in \calh,\\
\\
\dis w(0,x)=\phi_n(x),
\end{array}\right.
\end{equation}
where $\call$ is given by
\begin{equation}\label{eq:ell}
 \call[f](x)=\frac{1}{2} Tr \;GG^*\; \nabla^2f(x)
+ \< x,A^*\nabla f(x)\>.
\end{equation}
In particular $w_n$ satisfies the integral equation
\begin{equation}
  w_n(t,x) =R_{t}\phi_n](x)+\int_0^t\red{\left[ R_{t-s}[
H_{min}(\nabla^B w_n(s,\cdot)
](x)+\ell_{0}(s)\right]}\; ds,\qquad t\in [0,T].\
x\in \calh,\label{solmildHJB-forwardn}
\end{equation}
By Theorem \ref{teo su controllo feedback}
calling $v_n(t,x)=w_n(T-t,x)$ we have
\begin{equation}\label{identif-v_n}
  v_n(t,x) =V_n(t,x):=\inf_{u\in \calu} J_n(t,x;u).
\end{equation}
and there exists an optimal feedback control $u_n(s)=\psi_n(s,Y(s))$.
Moreover, by Lemma \ref{lm:approximation} we know that
$$
v_n(t,x) \overset{\calk}{\rightarrow} v(t,x).
$$
Now it is enough to prove that $V_n(t,x) \to V(t,x)$ pointwise.
Given $\eps >0$, we have, for $n$ large enough,
\begin{align*}
 V_n(t,x)&=\inf_{u\in\calu}
\left[\E \int_t^T \left(\red{\ell_{0}(s)}+\ell_1(u(s))\right)\;ds +\E  \phi_n(Y(T))\right]\\
&=\red{\inf_{u\in\calu}\left[\E \int_t^T \left(\red{\ell_{0}(s)}+\ell_1(u(s))\right)\;ds +\E \phi(Y(T))
+\E\left[\phi_n(Y(T)) -\phi(Y(T))\right]\right]}\\
&\leq\inf_{u\in\calu}\left[\E \int_t^T \left(\red{\ell_{0}(s)}+\ell_{0,n}(s,\cdot)+\ell_1(u(s))\right)\;ds +\E \phi(Y(T))\right]
+\varepsilon,
\end{align*}
where the last passage follows by the Dominated Convergence Theorem, and since $\phi$
and $\ell_{0,n}$ are uniformly continuous.
We have shown that
\[
 V_n(t,x)\leq V(t,x)+\varepsilon.
\]
Exchanging the role of $V_n$ and $V$ we also find that the reverse inequality holds true. Hence $V_n \to V$ pointwise an the claim follows.
\qed


\begin{remark}\label{rm:unif-cont-feedback}
 In Theorem \ref{teo:v=V} we have assumed further uniform continuity on the data.
When $U$ is compact the result still remain true if the data are only continuous.
\end{remark}


\end{document}